\documentclass[12pt]{article}
\def\C{\mathbf{C}}
\def\bC{\mathbf{\overline{C}}}
\def\CP{\mathbf{P}}
\title{On the spherical derivative of
a rational function} 
\author{Matthew Barrett and Alexandre Eremenko\thanks{Both authors
were supported by NSF grant DMS-1067886.}}
\begin{document}
\maketitle
\begin{abstract}
For a rational functions $f$ we consider the norm of the derivative 
with respect to the spherical metric and denote by $K(f)$ the
supremum of this norm. We give estimates of this quantity $K(f)$ both for 
an individual function and for sequences of iterates.


\end{abstract}

A rational function is a holomorphic map from the Riemann sphere
into itself. We equip the Riemann sphere with the usual spherical metric
whose length and area elements are
$$ds=\frac{|dz|}{1+|z|^2}\quad\mbox{and}\quad dA=\frac{dxdy}{(1+|z|^2)^2}.$$
So the norm of the derivative with respect to the spherical metric
is
$$\| f'\|(z):=|f'(z)|(1+|z|^2)/(1+|f(z)|^2).$$
In this paper we study the quantity 
$K(f)=\max_{\bC}\| f'\|$.
d'Ambra and Gromov \cite{d'ambra} proposed
to study the rate of growth of $\sup\|(f^n)'\|$
as $n\to \infty$ for the iterates $f^n$ of smooth maps
of Riemannian manifolds, especially those maps in a given class
for which this growth rate is the smallest possible.
Such maps are called ``slow''.
Slow maps of an interval and 
slow Hamiltonian diffeomorphisms of a $2$-torus
have been investigated in \cite{borichev2,PS} and \cite{borichev1}.

Let $f$ be a rational function of degree $d$.
As the map $f:\bC\to\bC$ is $d$-to $1$,
we conclude that
$$\int\int_{\C}\| f'\|^2 dA=
d(f)\int\int_\C dA.$$
This implies that
\begin{equation}\label{1}
K(f)\geq\sqrt{d(f)}
\end{equation}
We ask how small can $K(f)$ be for a function of given degree.
\vspace{.1in}

\noindent
1. It is known that $K(f)\geq 2$ for all rational functions of degree at
least $2$. In fact this holds for all smooth maps of the sphere
into itself which satisfy $\deg(f)\not\in\{0,1,-1\}$ \cite{Gromov}.
It is not known whether $K(f)=2$ can hold for rational functions
of degrees $3$ or $4$.  
\vspace{.1in}

\noindent
2. An interesting question is whether (\ref{1}) is best possible
in certain sense. We have
\vspace{.1in}

\noindent
{\bf Theorem 1.} {\em There exists an absolute constant $C$ with the following
property. For every $d\geq 2$ there exists a rational function of degree
exactly $d$ such that
$$K(f)\leq C\sqrt{d}.$$
}
\vspace{.1in}

An analogous result was obtained by Gromov \cite[Ch. 2D]{Gromov}
for smooth maps of spheres of arbitrary dimension.

Littlewood \cite{L} and Hayman \cite{Hayman}
studied the quantity
$$\phi(d)=\sup_{\deg f=d}\sup_{R>0}\frac{1}{R}\int\int_{|z|\leq R}
\frac{|f'(z)|}{1+|f(z)|^2}dxdy,$$
where the sup is taken over all rational functions of degree $d$.
For polynomials $\phi(d)$ was also studied in \cite{ES,E,LW}.
It is easy to see that $\phi(d)\leq\pi\sqrt{d}$ and Hayman obtained
$\phi(d)\geq c_1\sqrt{d}$ using a rational approximation of
elliptic functions. Our Theorem 1 implies this with a more elementary
proof.
Indeed, by a change of the independent variable,
$$\phi(d)=\sup_{\deg f=d}\int\int_{|z|\leq 1}\| f'\|\frac{dxdy}{1+|z|^2}.$$
Denote 
\begin{equation}\label{spherd}
f^\#=|f'|/(1+|f|^2).
\end{equation}
Let $f$ be the function from Theorem~1.
By rotating the sphere of the independent variable we may achieve that
$$\int\int_{|z|\leq 1}(f^\#)^2dxdy\geq \pi d/2,$$
because the spherical area of the image sphere is $\pi$ and it is covered
$d$ times. Let $M=\max_{|z|\leq 1}f^\#.$ Theorem~1 implies that
$M\leq C\sqrt{d}$,
so
\begin{eqnarray*}
\phi(d)&\geq&\int\int_{|z|\leq 1}f^\#dxdy\geq M^{-1}\int\int_{|z|\leq 1}
(f^\#)^2dxdy\\ \\
&\geq&C^{-1}d^{-1/2}\pi d/2=\pi C\sqrt{d}/2.
\end{eqnarray*}

{\em Proof of Theorem 1.}
For every positive integer $n$, consider the function 
$$f_n(z)=\prod_{k=-n}^n\tanh(z+2k).$$
We first show that
\begin{equation}\label{f}
\frac{|f_n|}{1+|f_n|^2}\leq C,
\end{equation}
where $C$ is independent of $n$.
We have
\begin{equation}\label{raz1}
||\tanh z|-1|\leq 4e^{-2x}<1,\quad x=|\Re z|\geq 1.
\end{equation}
Now fix $z\in\C$ and let $m$ be an integer such that
$|2m-\Re z|\leq 1$. For $|w-z|\leq 1/2$ put
$$g(w)=\left\{\begin{array}{ll}\tanh(w-2m),&\mbox{if}\quad|m|\leq n,\\
        1,&\mbox{otherwise}\end{array}\right.$$
Then evidently
\begin{equation}\label{dwa}
\frac{|g'(w)|}{1+|g(w)|^2}\leq\frac{|\tanh'(w)|}{1+|\tanh(w)|^2}\leq C_0.
\end{equation}
We write $f_n=gh$, $f_n^\prime=g'h+gh'$, and estimate $h$ using (\ref{raz1}):
$$|h(w)|\leq\prod_{k=0}^\infty\left(1+4e^{-2(2k+1/2)}\right)^2=:C_1.$$
Now $h$ is holomorphic in $|w-z|\leq 1/2$, so by Cauchy's theorem,
$$|h'(z)|\leq 4C_1.$$
Next we estimate $h(z)$ from below using (\ref{raz1}) again:
$$|h(z)|\geq \prod_{k=0}^\infty\left(1-4e^{-2(2k+1)}\right)^2=:C_2.$$
Combining all these estimates we obtain
$$\frac{|f^\prime_n(z)|}{1+|f_n(z)|^2}\leq
\frac{C_1|g'(z)|+4C_1|g(z)|}{1+C_2^2|g(z)|^2}
\leq C_2^{-2}(C_1C_0+4C_1).$$
This proves (\ref{f})
\vspace{.1in}

When $|\Re z|\geq n+1$, we will obtain better estimates. Using
$$|\tanh'(z)|=|\cosh(z)^2|\leq 16e^{-2|x|},\quad |x|=|\Re z|\geq 1,$$
we obtain for $z>2n$ and $\xi=\Re z-2n$:
$$|f^\prime_n(z-2n)|\leq C_1\sum_{k=0}^\infty|\cosh(\xi+2k)|^2\leq
16C_1\sum_{k=0}e^{-2(\xi+2k)}=16C_1C_3e^{-2\xi},$$
which gives
$$\frac{|^\prime_n(z)|}{1+|f_n(z)|^2}\leq Ce^{-2(|x|-n)}\leq Ce^{-2|x|/n}$$
if $|x|\geq n+2$. Combining this with (\ref{f}) we obtain
$$\frac{|f^\prime_n(z)|}{1+|f_n(z)|^2}\leq Ce^{-2|x|/n}$$
for all $z$.
As $f_n$ has period $\pi$, there
exists a rational function $R_n$
such that $f_n(nz)=R_n(e^{2z}).$
This rational function has degree $2n^2$ and the derivative
satisfies
$\| R_n^\prime\|\leq Cn.$
This completes the proof
of Theorem 1.
\vspace{.1in}

\noindent
{\bf Theorem 2.} {\em There exists an absolute constant $c>1$
with the property that
$$K(f)\geq c\sqrt{d}$$
for all rational functions of degree $d\geq 2$.}
\vspace{.1in}

This can be considered as an analog of a result of Tsukamoto \cite{Tsukamoto}.
He studied spherical derivatives  (\ref{spherd}) of
meromorphic functions $F:\Delta\to\bC$, where
$\Delta$ is the unit disc with the Euclidean metric and 
proved that there exists an absolute constant $c_1<1$
with the property that $\omega(F(\Delta))\leq c_1\pi$
for all meromorphic functions $F$
satisfying $F^\#\leq 1$,
where
$$\omega(F(\Delta))=\int_{|z|<1}(F^\#(z))^2dxdy.$$

We derive Theorem~2 from this result. In fact we show that Theorem~2
holds with $c=1/\sqrt{c_1}$.
\vspace{.1in}

{\em Proof of Theorem 2}. Proving by contradiction,
we suppose that there exists a sequence
$f_m$ of rational functions of degrees $m$,
such that 
\begin{equation}
\label{10}
K(f_m)/\sqrt{m}\to b<1/\sqrt{c_1}. 
\end{equation}
Let $\omega$ be the spherical area measure,
so that
$$\int_\bC d\omega=\pi,$$
and $\omega_m=f^*_m\omega$
the pull back of $\omega$ by $f_m$. Then
\begin{equation}\label{raz}
\int_\bC d\omega_m=\pi m.
\end{equation}
It is easy to see that we can find discs $D_m=D(a_m,r_m)\subset\bC$
(with respect to the spherical metric) of radii $r_m$
such that
\begin{equation}\label{13}
\int_{D_m}d\omega=\pi/(mb^2)\quad\mbox{and}\quad
\int_{D_m}d\omega_m\geq \pi/b^2.
\end{equation}
To show this, choose $r_m$ so that the first equation is satisfied
and then integrate 
$$F(a):=\int_{D(a,r_m)}d\omega_m(a)$$
with respect to $a$. Evidently $r_m\sim 1/b\sqrt{\pi m}, m\to\infty$.

Let $a_m^\prime$ be the point diametrically opposite to $a_m$, and
let $\phi_m:\C\to\bC\backslash\{ a_m^\prime\}$ be the conformal map
(inverse to a stereographic
projection) such that $\phi_m(0)=a,\phi_m(\Delta)=D_m$,
then 
\begin{equation}\label{11}
\phi_m^\#(z)\leq r_m(1+o(1)),\; m\to\infty,
\end{equation} uniformly
with respect to $z$.
Then $F_m=f_n\circ\phi_m$ is a normal family.
Let $F=\lim F_m$. From (\ref{10}), (\ref{11})
follows that then $F^\#\leq 1$, but the area of
$F(\Delta)\geq \pi/b^2$ in view of (\ref{13}),
contradicting the result of Tsukamoto.
\vspace{.1in}

Theorems 1 and 2 have analogs for maps $\CP^1\to\CP^n$ which are stated
and proved in the same way as for $n=1$, using he Fubini--Study metric
for the norm of the derivative. Constants $C$ and $c$ will of course
depend on $n$.
\vspace{.1in}

Now we consider dynamical questions. By $f^n$ we denote the $n$-th iterate,
and our standing assumption is that $d(f)\geq 2$.
We define
$$k_\infty(f)=\lim_{n\to\infty}\frac{1}{n}\log K(f^n).$$
The limit always exists because the sequence $a_n=\log K(f^n)$
is subadditive, $a_{m+n}\leq a_m+a_n$ and for every such positive sequence
the limit $\lim_{n\to\infty}a_n/n$ exists and is equal to $\inf_na_n/n$
(see, for example, Lemma 1.16 in \cite{Bowen}).

It follows that $k_\infty(f^m)=mk_\infty(f)$. 

Notice that $k_\infty$ is independent of the choice of a smooth
Riemannian metric on the sphere, and is invariant under
conjugation by conformal automorphisms.  
Obviously, $k_\infty(f)\leq \log K(f)$. 
\vspace{.1in}

\noindent
3. What is the smallest value of $k_\infty(f)$
for rational functions of given degree?
\vspace{.1in}

The trivial lower estimate of $K(f)$ gives
\begin{equation}\label{4}
k_\infty(f)\geq (1/2)\log d(f).
\end{equation}
We will see that
equality never happens, and that the Latt\'es functions are
not extremal for minimizing $k_\infty$.
For functions $f_d$ of degree $d$ from Theorem~1 we have
$$k_\infty(f_d)/\log d\to 1/2,\quad d\to\infty,$$
so the $(1/2)$ in (\ref{4}) cannot be replaced
with a larger constant.

In \cite{EL} these quantities were studied for polynomials,
in particular, the inequality $k_\infty(f)\geq \log d(f)$
was established for polynomials, with equality only if
$f$ is conjugate to $z^d$.

Let us consider a slightly different quantity,
the {\em maximum characteristic exponent}
\begin{equation}\label{chim}
\chi_m(f)=\sup_z\limsup_{n\to\infty}
\frac{1}{n}\log\|(f^n)'(z)\|\leq k_\infty(f).
\end{equation}
The difference in the definitions of $k_\infty$ and $\chi_m$ is
in the order of $\max_z$ and $\lim_{n\to\infty}$.
Przytycki proved in 1998 (reproduced in \cite{Pr}) that the same quantity
$\chi_m(f)$ can be obtained by taking the sup over periodic pints $z$ of
$f$, in which case the $\limsup$ in (\ref{chim}) can be of course replaced
by the ordinary limit. Moreover, he proved
the following:
\vspace{.1in}

\noindent
{\bf Theorem P.} {\em For every $\epsilon>0$ there exists a periodic point
$z$ such that 
\begin{equation}
\label{12}
\frac{1}{m}\log\|(f^m)'(z)\|\geq k_\infty-\epsilon,
\end{equation}
where $m$ is a period of $z$.

In particular, $k_\infty=\chi_m$, and one can replace $\sup_{z\in\bC}$ in (\ref{chim})
by $\sup$ over all periodic points.}
\vspace{.1in}

4. Let $\chi_a(f)$ be the average value of
$$\log\|(f^n)'(z)\|$$
over the measure $\mu$ of maximal entropy.
According to the multiplicative ergodic theorem, we have
$$\chi_a(f)=\lim_{n\to\infty}\frac{1}{n}\log\|(f^n)'(z)\|$$
almost everywhere with respect to $\mu$.

\vspace{.1in}

\noindent
{\bf Theorem 4.} {\em $\chi_a(f)\geq (1/2)\log d$,
this is best possible and equality holds only for Latt\'es
examples.}
\vspace{.1in}

{\em Proof.}
The estimate follows from the formula
$$\chi_a(f)=h(f)/\dim\mu,$$
where $h(f)=\log d$ is the topological entropy, and
$\dim\mu$ is the Hausdorff dimension of the maximal
measure, see, for example \cite{ELy} and references therein.
Obviously
$\dim\mu\leq 2$ so we obtain our inequality.
On the other hand, a theorem of A. Zdunik \cite{Zdunik} says that
$\dim\mu=2$ can happen only for Latt\'es examples. This completes the
proof of the theorem.
\vspace{.1in}

Now we consider the Latt\'es functions \cite{Milnor}.
\vspace{.1in}

\noindent
{\bf Proposition 1.} {\em If $L$ is a Latt\'es function then
$\chi_m(L)\geq \log d(L)$.}
\vspace{.1in}

{\em Proof.} A Latt\'es function can be defined by a functional equation
$$F(\lambda z)=L\circ F(z),\quad d(L)=\lambda^2,$$
where $F$ is an elliptic function with a critical point at $0$.
Assuming without loss of generality that $F(0)=0$ we conclude 
that $0$ is a fixed point of $L$. Now the derivative at this fixed
point is
$$\lambda\lim_{z\to 0}F'(\lambda z)/F'(z)=\lambda^2=d(L),$$
from which follows that $\chi_m(L)\geq d(L)$.
\vspace{.1in}

To summarize, we have 4 quantities satisfying inequalities
\begin{equation}\label{5}
\frac{1}{2}\log d\leq\chi_a(f)\leq\chi_m(f)=k_\infty(f)
\leq \log K(f).
\end{equation}
For Latte\'es functions
we have
$$(1/2)\log d=\chi_a(L)<\chi_m(L)=\log d.$$
In the first inequality equality holds only for Latt\'es functions.
For functions $f_d$ constructed in Theorem~1, we have 
$\log K(f_d)\leq (1/2)\log d+O(1).$

Which of the rest of these inequalities (\ref{5}) are strict and what
are the conditions of equality? According to a private communication
of Przytycki, the method of Zdunik can be used to show that
$\chi_a(f)=\chi_m(f)$ is only possible when $f$ is conjugate to
$z^{\pm d}$ in which case both quantities are equal to $\log d$. 
\vspace{.1in}

\noindent
{\bf Proposition.} {\em Let $M=\{ z:\| f'(z)\|=K(f)\}$.
Then $k_\infty=\log K(f)$ if and only if $M$
contains a cycle.
}
\vspace{.1in}

{\em Proof.} If $z$ is a point whose trajectory
is in $M$, then
$\| (f^m)'(z)\|=\| f'(z)\|^m$ so $k_\infty(f)=\log K(f)$.

Suppose now that $k_\infty=\log K(f)$. We claim that 
$$\bigcap_{j=0}^\infty f^j(M)\neq\emptyset.$$
Indeed, otherwise for some $m$ we have
$$\bigcap_{j=0}^m f^j(M)=\emptyset,$$
and if this holds, then
\begin{eqnarray*}
mk_\infty(f)&=&k_\infty(f^m)\leq \log K(f^m)=\log\max_z
\prod_{j=0}^{m-1}\| f'(f^j(z)\|\\
&<& m\log K(f),
\end{eqnarray*}
contrary to our assumption. This proves the claim.

The set $M_{\infty}$ is forward invariant.
All sets
$$\bigcap_{j=0}^m f^j(M)$$
are real algebraic subsets of $\bC$, and thus the set $M_\infty$
is also real algebraic. If $M_\infty$ is finite it must
contain a cycle. 

Suppose that $M$ is infinite. Let $M_1,...M_q$ be the list
of connected components of $M$. Then there must be a cycle of components.
Let $M_1,\ldots,M_k$ be this cycle, so that $f(M_j)\subset M_{j+1}$, 
$f(M_k)=M_1$.
If any of these components contains singular points, then all singular
points in these components form a finite invariant set, and this set must
contain a cycle. If all $M_j,$ $1\leq j\leq k$ are smooth,  
then $f^k:M_1\to M_1$ is an expanding map of a circle, so it must have
a fixed point. Again we obtain a cycle in $M$.
This completes the proof.

\vspace{.1in}

{\em Purdue University

West Lafayette, IN 47907

eremenko{@}math.purdue.edu}
\end{document}